\documentclass[11pt]{amsart}
\numberwithin{equation}{section}
\newtheorem{theorem}{Theorem}

\newtheorem{corollary}{Corollary}
\numberwithin{theorem}{section} \numberwithin{lemma}{section}
 \numberwithin{definition}{section}
\numberwithin{proposition}{section}
\newtheorem{remark}{Remark}[section]
\def\R{\bf R}

\def\al{\aligned}
\def\eal{\endaligned}
\def\M{{\bf M}}
\def\be{\begin{equation}}
\def\ee{\end{equation}}
\def\lab{\label}

\def\R{\bf R}
\def\M{\bf M}

\def\al{\aligned}
\numberwithin{equation}{section}

\begin{document}

\tracingpages 1
\title[ancient solutions]{\bf A note on time analyticity for ancient solutions of the heat equation }
\author{Qi S. Zhang}
\address{ Department of
Mathematics, University of California, Riverside, CA 92521, USA }
\date{May 2019}

\begin{abstract}
It is well known that generic solutions of the heat equation are not analytic in time in general. Here it is proven that ancient solutions  with exponential growth are analytic in time in ${\M} \times (-\infty, 0]$. Here  $\M=\R^n$ or is a manifold with Ricci curvature bounded from below. Consequently  a necessary and sufficient condition is found on the solvability of backward heat equation in the class of functions with exponential growth.
\end{abstract}
\maketitle
%\tableofcontents
\section{Statement of result and proof}

The goal of the note is to prove time analyticity of certain generic ancient solutions to the heat equation on ${\M} \times (-\infty, 0]$ where $\M=\R^n$ or some noncompact Riemannian manifolds. This is somewhat unexpected since it is well known that generic solutions to the forward heat equation are not analytic in time. Even if $\M=\R^1$ and $u$ is analytic in $x$ when $t>0$, time analyticity requires at least the initial value is analytic in $x$, c.f. \cite{Wi:1} Corollary 3.1 b. Ancient solutions of evolution equations are important since they not only represent structures of solutions near a high value but also can be regarded as solutions of the backward equation.
 Backward heat equations play important roles in stochastic analysis and Ricci flows e.g. It is known that the Cauchy problem to the backward heat equation is not solvable in general. One application of the main result is a necessary and sufficient condition for solvability if the solutions grow no faster than exponential functions.
 It is expected that the phenomenon observed here can be extended to many other evolution equations.

Let us recall some relevant results for ancient solutions of the heat equation. Let $\M=\R^n$ or complete noncompact Riemannian manifold with nonnegative curvature.
It is proven in \cite{SZ:1} that sublinear ancient solutions are constants.  In \cite{LZ:1} it is proven that the space of  ancient solutions of polynomial growth has finite dimension and the solutions are polynomials in time. In the paper \cite{CM:1}, a sharp dimension estimate of the space is given. See the papers \cite{Ca:1} and \cite{CM:2} for applications to the study of mean curvature flows, and also \cite{Ha:1} in the graph case. Earlier, backward heat equations have been studied by many authors, see \cite{Mi:1}, \cite{Yo:1} e.g, and treated in many text books. A necessary and sufficient solvability criteria seems lacking, except when 
$\M$ is a bounded domain for which semigroup theory gives an abstract criteria \cite{CJ:1} Theorem 9.

In order to state the result, let us first introduce a bit of notations.
We use $\M$ to denote a $n$ dimensional, complete, noncompact Riemannian manifold, $Ric$ to denote the Ricci curvature and $0$ a reference point on $\M$, $d(x, y)$ is the geodesic distance for $x, y \in \M$.

\begin{theorem}
\lab{thhepoly}
 Let $\M$ be a complete, n dimensional, noncompact Riemannian manifold such that the Ricci curvature satisfies $Ric \ge -(n-1) K_0$ for a nonnegative constant $K_0$.

 Let $u$ be a smooth, ancient solution of the heat equation $\Delta u - \partial_t u=0$ of exponential growth, namely
\be
\lab{expg}
 |u(x, t)| \le A_1 e^{ A_2 ( d(x, 0) + |t|)}, \quad \forall (x, t) \in {\M} \times (-\infty, 0],
\ee where $A_1, A_2$ are positive constants and $0$ is a reference point on $\M$.
Then $u=u(x, t)$ is analytic in $t$ with radius $\infty$, moreover,
\be
\lab{utaylor}
u(x, t)= \sum^\infty_{j=0} a_j(x) \frac{ t^j}{j !},
\ee with $\Delta a_{j}(x) = a_{j+1}(x)$. In addition,
\be
\lab{thaijie}
|a_{j}(x)| \le
A_3 e^{ A_4 (j + d(x, 0))}, \quad j=0, 1,2,...;
\ee where  $A_3, A_4$ are positive constants depending on $A_1, A_2, K_0, n$, and $A_3$ also depends  on $|B(0, 1)|^{-1}$.
\proof
\end{theorem}

Let us start with a well known parabolic mean value inequality which can be found in Theorem 14.7 of \cite{Li:1} e.g.  Suppose $u$ is a positive subsolution to the heat equation on ${\M} \times [0, T]$. Let $T_1, T_2 \in [0, T]$ with $T_1<T_2$, $R>0$, $p>0$, $\delta, \eta \in (0, 1)$. Then there exist  positive constants $C_1$ and $C_2$, depending only on $p, n$ such that
\be
\lab{pmvi}
\al
\sup_{B(0, (1-\delta) R) \times [T_1, T_2]} u^p &\le C_1 \frac{\bar{V}(2R)}{|B(0, R)|}
( R \sqrt{K_0} + 1) \exp ( C_2 \sqrt{K_0 (T_2-T_1)}) \\
&\qquad \times
\left(\frac{1}{\delta R} + \frac{1}{\eta T_1} \right)^{n+2} \,
\int^{T_2}_{(1-\eta) T_1} \int_{B(0, R)} u^p(y, s) dy ds;
\eal
\ee Here $\bar{V}(R)$ is the volume of geodesic balls of radius $R$ in the simply connected space form with constant sectional curvature $-K_0$; $|B(0, R)|$ is the volume of the geodesic ball $B(0, R)$ with center $0$ and radius $R$.

Let $u$ be an ancient solution to the heat equation. Then $u^2$ is a subsolution.  Given a
positive integer $j$, by shifting the time suitably, we can apply the mean value inequality to deduce
\be
\lab{mviq}
\sup_{B(0, j) \times [-j, 0]} u^2 \le C_1 e^{C_2 \sqrt{K_0} j}
\int^{0}_{-(j+1)} \int_{B(0, j+1)} u^2(y, s) dy ds;
\ee where we have used the Bishop-Gromov volume comparison theorem. Note that the constants $C_1$, $C_2$ may have changed and $C_1$ now also depends on the lower bound of $|B(0, 1)|$, i.e. the volume noncollapsing constant.  For each $j$, observe that $\partial^j_t u$ is also a solution of the heat equation. Applying (\ref{mviq}) with $u$ replaced by  $\partial^j_t u$, we deduce
\be
\lab{mviqj}
\sup_{Q(0, j)} (\partial^j_t u)^2 \le C_1 e^{C_2 \sqrt{K_0} j} \int_{Q(0, j+1)} (\partial^j_s u)^2(y, s) dy ds,
\ee where we have used the notation $Q(0, j) \equiv B(0, j) \times [-j, 0]$ to denote the space time cube of size $j$ with vertex $(0, 0)$. Note that this is not the standard parabolic cube since spatial and time scale is the same.

  Denote by $\psi$ a standard Lipschitz cut off function supported in $Q(0, j+0.5)$ such that $\psi=1$ in $Q(0, j)$ and $|\nabla \phi|^2 + |\partial_t \psi| \le C$. Since $u$ is a smooth solution to the heat equation, we compute
\[
\al
\int_{Q(0, j+0.5)} &(\Delta u)^2 \psi^2 dxdt = \int_{Q(0, j+0.5)} u_t \Delta u \psi^2 dxdt\\
&=-\int_{Q(0, j+0.5)} ((\nabla u)_t \nabla u) \, \psi^2 dxdt - \int_{Q(0, j+0.5)} u_t \nabla u \nabla \psi^2 dxdt\\
&= - \frac{1}{2} \int_{Q(0, j+0.5)} (|\nabla u |^2)_t \, \psi^2 dxdt - 2 \int_{Q(0, j+0.5)} u_t \psi \nabla u \nabla \psi dxdt\\
&\le \frac{1}{2} \int_{Q(0, j+0.5)} |\nabla u |^2 \, (\psi^2)_t dxdt +\frac{1}{2} \int_{Q(0, j+0.5)} (u_t)^2 \psi^2 dxdt \\
&\qquad + 2 \int_{Q(0, j+0.5)} |\nabla u|^2 |\nabla \psi|^2 dxdt.
\eal
\]Therefore
\[
\int_{Q(0, j+0.5)} (\Delta u)^2 \psi^2 dxdt \le C \int_{Q(0, j+0.5)} |\nabla u |^2 dxdt.
\]This and the standard Cacciopoli inequality (energy estimate)  between the cubes $Q(0, j+0.5)$ and  $Q(0, j+1)$ show that
\be
\lab{ddu2r-4}
\int_{Q(0, j)} (\Delta u)^2  dxdt \le C_0 \int_{Q(0, j+1)} u^2  dxdt.
\ee Here $C_0$ is a universal constant.

Since $\partial^j_t u$ is also an ancient solution, we can replace $u$ in (\ref{ddu2r-4}) by $\partial^j_t u$ to deduce
\[
\int_{Q(0, j+1)} (\partial^j_t u)^2  dxdt = \int_{Q(0, j+1)} (\Delta  \partial^{j-1} u)^2  dxdt\\
\le C_0 \int_{Q(0, j+2)} (\partial^{j-1}_t u)^2  dxdt
\]By induction, we deduce
\be
\lab{utk<}
\int_{Q(0, j+1)} (\partial^j_t u)^2  dxdt \le C^j_0 \int_{Q(0, 2 j+1)} u^2  dxdt.
\ee

Substituting (\ref{utk<}) to (\ref{mviqj}), we
find that
\[
\sup_{Q(0, j)} (\partial^i_t u)^2 \le C_1 e^{C_2 \sqrt{K_0} j} C^j_0 \int_{Q(0, 2j+1)} u^2  dxdt,
\] for all $i=0, 1, 2, ... j.$  This implies, by the exponential growth condition (\ref{expg}) and the Bishop-Gromov volume comparison theorem, that
\be
\lab{djtu}
\sup_{Q(0, j)} |\partial^i_t u| \le C_1 C^j_3
\ee for all $i=0, 1, 2, ... j.$ Here $C_3$ is a positive constant.

Fixing a number $R \ge 1$, for $x \in B(0, R)$, choose an integer $j \ge R$ and $t \in [- j^2, 0]$. Taylor's theorem implies that
\be
\lab{jtaylor}
u(x, t)- \sum^{j-1}_{i=0} \partial^i_t u(x, 0) \frac{ t^i}{i !} = \frac{t^j}{ j !} \partial^j_s u(x, s),
\ee where $s=s(x, t, j) \in [t, 0]$.
By (\ref{djtu}), the right hand side of (\ref{jtaylor}) converges to $0$ uniformly on $Q(0, R)$
as $j \to \infty$. Hence
\be
\lab{ujiexi}
u(x, t)= \sum^{\infty}_{i=0} \partial^j_t u(x, 0) \frac{ t^j}{j !}
\ee i.e. $u$ is analytic in $t$ with radius $\infty$.  Writing $a_j= a_j(x) = \partial^j_t u(x, 0)$. By (\ref{djtu}) again, we have
\be
\partial_t u(x, t) = \sum^{\infty}_{j=0} a_{j+1}(x) \frac{ t^j}{i !},
\ee
\be
\Delta u(x, t) = \sum^{\infty}_{j=0} \Delta a_j(x) \frac{ t^j}{j !},
\ee where both series converge uniformly on $Q(0, R)$ for any fixed $R>0$.
Since $u$ is a solution of the heat equation,
this implies
\[
\Delta a_j(x) =a_{j+1}(x)
\] with
\[
|a_j(x)| \le A_3 e^{ A_4 (j + d(x, 0))}.
\] Here $A_3$ and $A_4$ are positive constants with $A_3$ depending on $|B(0, 1)|^{-1}$. This completes the proof of the theorem. \qed
\medskip

An immediate application is the following:

\begin{corollary}
\lab{cobhe}
Let $\M$ be as in the theorem. Then the Cauchy problem for the backward heat equation
\be
\lab{bhe}
\begin{cases} \Delta u + \partial_t u = 0, \quad {\M} \times[0, \infty);\\
u(x, 0)=a(x)
\end{cases}
\ee has a smooth  solution of exponential growth if and only if
\be
\lab{aijie}
|\Delta^j a(x)| \le A_3 e^{ A_4 (j + d(x, 0))}, \quad j=0, 1,2,...;
\ee where  $A_3, A_4$ are positive constants depending on $A_1, A_2, K_0, n$, and $A_3$ also depends  on $|B(0, 1)|^{-1}$.
\proof
\end{corollary}

Suppose (\ref{bhe}) has a smooth solution of exponential growth, say $u=u(x, t)$. Then
$u(x, -t)$ is an ancient solution of the heat equation with exponential growth. By the theorem
\[
u(x, -t) = \sum^\infty_{j=0} \Delta^j a(x) \frac{(-t)^j}{j!}
\] Then (\ref{aijie}) follows from the theorem since $\Delta^j a(x) = a_j(x)$ in the theorem.

On the other hand, suppose (\ref{aijie}) holds. Then it is easy to check that
\[
u(x, t) = \sum^\infty_{j=0} \Delta^j a(x) \frac{t^j}{j!}
\]is a smooth ancient solution of the heat equation. Indeed, the bounds (\ref{aijie}) guaranty that the above series and the series
\[
 \sum^\infty_{j=0} \Delta^{j+1} a(x) \frac{t^j}{j!}, \quad
\sum^\infty_{j=0} \Delta^{j} a(x)  \frac{\partial_t t^j}{j!},
\]all converges absolutely and uniformly in $Q(0, R)$ for any fixed $R>0$. Hence $\Delta u - \partial_t u =0$.
 Moreover $u$ has exponential growth since
 \[
|u(x, t)| \le \sum^\infty_{j=0} |\Delta^j a(x)| \frac{|t|^j}{j!} \le A_3  e^{A_4 d(x, 0)} \sum^\infty_{j=0} \frac{\left(e^{A_4} |t| \right)^j}{j!} = A_3  e^{A_4 d(x, 0) + e^{A_4} |t|}.
\]
Thus $u(x, -t)$ is a solution to the Cauchy problem of the backward heat equation (\ref{bhe}) of exponential growth.
\qed

{\remark  For the conclusion of the theorem to hold, some growth condition for the solution is necessary. Tychonov's non-uniqueness example can be modified as follows. Let $v=v(x, t)$ be Tychonov's solution of the heat equation in ${\R^n} \times (-\infty, \infty)$,  which is $0$ when $t \le 0$ but nontrivial for $t>0$.
Then $u \equiv v(x, t+1)$ is a nontrivial ancient solution. It is clearly not analytic in time.
Note that $|u(x, t)|$ grows faster than $e^{c |x|^2}$ for some $x$ and $t$. We are not sure if this bound is sufficient for time analyticity.}

{\remark  If $\M =\R^n$, then the ancient solution in the theorem is also analytic in space variables. For general manifolds, space analyticity requires certain bounds on curvature and its derivatives.}

\medskip
{\bf Acknowledgment.} We wish to thank Professors Bobo Hua and F. H. Lin for helpful discussions. Part of the paper was written when the author was visiting the School of Mathematics at Fudan University. We are grateful to Professor Lei Zhen the invitation and warm hospitality, to the Simons foundation for its support.

\bigskip

\noindent e-mails:
qizhang@math.ucr.edu


\begin{thebibliography}{00}



\bibitem[Ca1]{Ca:1}  M. Calle, {\it Bounding dimension  of ambient  space by density for mean  curvature flow},  Math. Z. 252(2006), no. 3, 655-668.




\bibitem[CJ]{CJ:1}
Christensen, Ann-Eva;  Johnsen, Jon {\it  Final value problems for parabolic differential equations and their well-posedness}, arXiv:1707.02136

\bibitem [CM1]{CM:1} Colding, Tobias H.; Minicozzi, William P., II
{\it Optimal bounds for ancient caloric functions}, arXiv:1902.01736


\bibitem [CM2]{CM:2} Colding, Tobias H.; Minicozzi, William P., II
{\it Complexity of parabolic systems},  arXiv:1903.03499


%\bibitem [Gr]{Gr:1}Grigoryan, A. A. {\it The heat equation on noncompact Riemannian manifolds.} (Russian) Mat. Sb. 182 (1991), no. 1, 55-87; translation in Math. USSR-Sb. 72 (1992), no. 1, 47-77

\bibitem [Ha]{Ha:1} Hua, Bobo, {\it
Dimensional bounds for ancient caloric functions on graphs }, arXiv:1903.02411.







\bibitem [Li] {Li:1} Li, Peter,
{\it Geometric analysis.} Cambridge Studies in Advanced Mathematics, 134. Cambridge University Press, Cambridge, 2012. x+406 pp.

%\bibitem [Ln]{Ln:1} Lin, Fang-Hua, {\it A uniqueness theorem for parabolic equations.} Comm. Pure Appl. Math. 43 (1990), no. 1, 127-136


\bibitem [LZ]{LZ:1}  Lin, Fanghua and Zhang, Qi S. {\it  On ancient solutions of the heat equation},
arXiv:1712.04091, Comm. Pure Applied Math. to appear 2019.


\bibitem [Mi]{Mi:1} Miranker, W. L. {\it A well posed problem for the backward heat equation}. Proc. Amer. Math. Soc. 12 1961 243-247.



\bibitem[SZ]{SZ:1} Souplet, Philippe; Zhang, Qi S. {\it Sharp gradient estimate and Yau's Liouville theorem for the heat equation on noncompact manifolds.} Bull. London Math. Soc. 38 (2006), no. 6, 1045-1053.



\bibitem [Yo]{Yo:1} Yosida, Kosaku, {\it  An abstract analyticity in time for solutions of a diffusion equation},     Proc. Japan Acad.
    Volume 35, Number 3 (1959), 109-113.

\bibitem [Wi]{Wi:1}Widder, D. V.
{\it Analytic solutions of the heat equation.}
Duke Math. J. 29 (1962) 497-503.



\end{thebibliography}
\enddocument